\input amstex
\let\myfrac=\frac
\input eplain
\let\frac=\myfrac
\input epsf




\loadeufm \loadmsam \loadmsbm
\message{symbol names}\UseAMSsymbols\message{,}

\font\myfontdefault=cmr10

\font\mytdmchapfont=cmb10 at 14pt
\font\mytdmheadfont=cmb10 at 10pt
\font\mytdmsubheadfont=cmr10
	
\magnification 1200
\newif\ifinappendices
\newif\ifundefinedreferences
\newif\ifchangedreferences
\newif\ifloadreferences
\newif\ifmakebiblio
\newif\ifmaketdm

\undefinedreferencesfalse
\changedreferencesfalse


\loadreferencestrue
\makebibliofalse
\maketdmfalse

\def\headpenalty{-400}     
\def\proclaimpenalty{-200} 

%
%

\def\alphanum#1{\ifcase #1 _\or A\or B\or C\or D\or E\or F\or G\or H\or I\or J\or K\or L\or M\or N\or O\or P\or Q\or R\or S\or T\or U\or V\or W\or X\or Y\or Z\fi}
\def\gobbleeight#1#2#3#4#5#6#7#8{}

\newwrite\references
\newwrite\tdm
\newwrite\biblio

\newcount\chapno
\newcount\headno
\newcount\subheadno
\newcount\procno
\newcount\figno
\newcount\citationno

\def\setcatcodes{%
\catcode`\!=0 \catcode`\\=11}%

\ifloadreferences
    {\catcode`\@=11 \catcode`\_=11%
    \global\def\_@citation@LabB{1}
\global\def\_@citation@LabA{2}
\global\def\_@citation@RosSpruck{3}
\global\def\_@citation@SmiB{4}
\global\def\_@citation@SmiC{5}
\global\def\_@citation@SmiG{6}
\global\def\_@proc@Smooth_Foliation{1.1}
\global\def\_@proc@First_Derivative_Of_A{2.1}
\global\def\_@proc@Diff_Eq_For_Equidistant_Family{2.2}
\global\def\_@proc@Growth_Is_Assured{2.3}
\global\def\_@proc@Curvature_Of_Equidistant_Surfaces_Bounded_Above{2.4}
\global\def\_@proc@It_Is_A_Foliation_Of_Something{3.1}
\global\def\_@proc@The_Foliation_Is_Complete_At_Infinity{4.1}
\global\def\_@proc@Weak_Lower_Bound_Of_Curv_For_Level_Surfaces{4.3}
\global\def\_@proc@The_Foliation_Is_Complete_At_Zero{4.4}
\global\def\_@proc@The_Foliation_Is_Smooth{5.1}
    }%
\else
    \openout\references=references.tex
\fi

\newcount\newchapflag 
\newcount\showpagenumflag 

\global\chapno = -1 
\global\citationno=0
\global\headno = 0
\global\subheadno = 0
\global\procno = 0
\global\figno = 0

\def\resetcounters{%
\global\headno = 0%
\global\subheadno = 0%
\global\procno = 0%
\global\figno = 0%
}

\global\newchapflag=0 
\global\showpagenumflag=0 

\def\chinfo{\ifinappendices\alphanum\chapno\else\the\chapno\fi}%
\def\headinfo{\ifinappendices\alphanum\headno\else\the\headno\fi}%
\def\subheadinfo{\headinfo.\the\subheadno}
\def\procinfo{\headinfo.\the\procno}
\def\figinfo{\headinfo.\the\figno}
\def\citationinfo{\the\citationno}%
\def\nextheadno{\global\advance\headno by 1 \global\subheadno = 0 \global\procno = 0}
\def\nextsubheadno{\global\advance\subheadno by 1}
\def\nextprocno{\global\advance\procno by 1 \procinfo}
\def\nextfigno{\global\advance\figno by 1 \figinfo}

{\global\let\noe=\noexpand%
%
%
\catcode`\@=11%
\catcode`\_=11%
\setcatcodes%
!global!def!_@@internal@@makeref#1{%
!global!expandafter!def!csname #1ref!endcsname##1{%
!csname _@#1@##1!endcsname%
!expandafter!ifx!csname _@#1@##1!endcsname!relax%
    !write16{#1 ##1 not defined - run saving references}%
    !undefinedreferencestrue%
!fi}}%
!global!def!_@@internal@@makelabel#1{%
!global!expandafter!def!csname #1label!endcsname##1{%
!edef!temptoken{!csname #1info!endcsname}%
!ifloadreferences%
    !expandafter!ifx!csname _@#1@##1!endcsname!relax%
        !write16{#1 ##1 not hitherto defined - rerun saving references}%
        !changedreferencestrue%
    !else%
        !expandafter!ifx!csname _@#1@##1!endcsname!temptoken%
        !else
            !write16{#1 ##1 reference has changed - rerun saving references}%
            !changedreferencestrue%
        !fi%
    !fi%
!else%
    !expandafter!edef!csname _@#1@##1!endcsname{!temptoken}%
    !edef!textoutput{!write!references{\global\def\_@#1@##1{!temptoken}}}%
    !textoutput%
!fi}}%
!global!def!makecounter#1{!_@@internal@@makelabel{#1}!_@@internal@@makeref{#1}}%
!unsetcatcodes%
}
\makecounter{ch}%
\makecounter{head}%
\makecounter{subhead}%
\makecounter{proc}%
\makecounter{fig}%
\makecounter{citation}%
\def\newref#1#2{%
\def\temptext{#2}%
\edef\bibliotextoutput{\expandafter\gobbleeight\meaning\temptext}%
\global\advance\citationno by 1\citationlabel{#1}%
\ifmakebiblio%
    \edef\fileoutput{\write\biblio{\noindent\hbox to 0pt{\hss$[\the\citationno]$}\hskip 0.2em\bibliotextoutput\medskip}}%
    \fileoutput%
\fi}%
\def\cite#1{%
$[\citationref{#1}]$%
\ifmakebiblio%
    \edef\fileoutput{\write\biblio{#1}}%
    \fileoutput%
\fi%
}%
%
%
%

\let\mypar=\par


\def\raggedleft{\leftskip=0pt plus 1fil \parfillskip=0pt}


\font\lettrinefont=cmr10 at 28pt
\def\lettrine #1[#2][#3]#4%
{\hangafter -#1 \hangindent #2
\noindent\hskip -#2 \vtop to 0pt{
\kern #3 \hbox to #2 {\lettrinefont #4\hss}\vss}}

\font\mylettrinefont=cmr10 at 28pt
\def\mylettrine #1[#2][#3][#4]#5%
{\hangafter -#1 \hangindent #2
\noindent\hskip -#2 \vtop to 0pt{
\kern #3 \hbox to #2 {\mylettrinefont #5\hss}\vss}}


\edef\Pagetitle={Blank}

\headline={\hfil\Pagetitle\hfil}

\footline={{\hfil\myfontdefault\folio\hfil}}

\def\nextoddpage
{
\newpage%
\ifodd\pageno%
\else%
    \global\showpagenumflag = 0%
    \null%
    \vfil%
    \eject%
    \global\showpagenumflag = 1%
\fi%
}


\def\newchap#1#2%
{%
%
%
\global\advance\chapno by 1%
\resetcounters%
%
%
\newpage%
\ifodd\pageno%
\else%
    \global\showpagenumflag = 0%
    \null%
    \vfil%
    \eject%
    \global\showpagenumflag = 1%
\fi%
\global\newchapflag = 1%
\global\showpagenumflag = 1%
%
%
{\font\chapfontA=cmsl10 at 30pt%
\font\chapfontB=cmsl10 at 25pt%
\null\vskip 5cm%
{\chapfontA\raggedleft\hfil%
{%
\ifnum\chapno=0
    \phantom{%
    \ifinappendices%
        Annexe \alphanum\chapno%
    \else%
        \the\chapno%
    \fi}%
\else%
    \ifinappendices%
        Annexe \alphanum\chapno%
    \else%
        \the\chapno%
    \fi%
\fi%
}%
\par}%
\vskip 2cm%
{\chapfontB\raggedleft%
\lineskiplimit=0pt%
\lineskip=0.8ex%
\hfil #1\par}%
\vskip 2cm%
}%
\edef\Pagetitle{#2}%
%
%
\ifmaketdm%
    \def\temp{#2}%
    \def\tempbis{\nobreak}%
    \edef\chaptitle{\expandafter\gobbleeight\meaning\temp}%
    \edef\mynobreak{\expandafter\gobbleeight\meaning\tempbis}%
    \edef\textoutput{\write\tdm{\bigskip{\noexpand\mytdmchapfont\noindent\chinfo\ - \chaptitle\hfill\noexpand\folio}\par\mynobreak}}%
\fi%
\textoutput%
}


\def\newhead#1%
{%
\ifhmode%
    \mypar%
\fi%
\ifnum\headno=0%
\ifinappendices
    \nobreak\vskip -\lastskip%
    \nobreak\vskip .5cm%
\fi
\else%
    \nobreak\vskip -\lastskip%
    \nobreak\vskip .5cm%
\fi%
\nextheadno%
\ifmaketdm%
    \def\temp{#1}%
    \edef\sectiontitle{\expandafter\gobbleeight\meaning\temp}%
    \edef\textoutput{\write\tdm{\noindent{\noexpand\mytdmheadfont\quad\headinfo\ - \sectiontitle\hfill\noexpand\folio}\par}}%
    \textoutput%
\fi%
\font\headfontA=cmbx10 at 14pt%
{\headfontA\noindent\headinfo\ -\ #1.\hfil}%
\nobreak\vskip .5cm%
}%


\def\newsubhead#1%
{%
\ifhmode%
    \mypar%
\fi%
\ifnum\subheadno=0%
\else%
    \penalty\headpenalty\vskip .4cm%
\fi%
\nextsubheadno%
\ifmaketdm%
    \def\temp{#1}%
    \edef\subsectiontitle{\expandafter\gobbleeight\meaning\temp}%
    \edef\textoutput{\write\tdm{\noindent{\noexpand\mytdmsubheadfont\quad\quad\subheadinfo\ - \subsectiontitle\hfill\noexpand\folio}\par}}%
    \textoutput%
\fi%
\font\subheadfontA=cmsl10 at 12pt
{\subheadfontA\noindent\subheadinfo\ #1.\hfil}%
\nobreak\vskip .25cm %
}%

%
%


\font\mathromanten=cmr10
\font\mathromanseven=cmr7
\font\mathromanfive=cmr5
\newfam\mathromanfam
\textfont\mathromanfam=\mathromanten
\scriptfont\mathromanfam=\mathromanseven
\scriptscriptfont\mathromanfam=\mathromanfive
\def\mathroman{\fam\mathromanfam}


\font\sansseriften=cmss10
\font\sansserifseven=cmss7
\font\sansseriffive=cmss5
\newfam\sansseriffam
\textfont\sansseriffam=\sansseriften
\scriptfont\sansseriffam=\sansserifseven
\scriptscriptfont\sansseriffam=\sansseriffive
\def\mathsf{\fam\sansseriffam}


\font\boldten=cmb10
\font\boldseven=cmb7
\font\boldfive=cmb5
\newfam\mathboldfam
\textfont\mathboldfam=\boldten
\scriptfont\mathboldfam=\boldseven
\scriptscriptfont\mathboldfam=\boldfive
\def\mathbf{\fam\mathboldfam}


\font\mycmmiten=cmmi10
\font\mycmmiseven=cmmi7
\font\mycmmifive=cmmi5
\newfam\mycmmifam
\textfont\mycmmifam=\mycmmiten
\scriptfont\mycmmifam=\mycmmiseven
\scriptscriptfont\mycmmifam=\mycmmifive

\def\hexa#1{\ifcase #1 0\or 1\or 2\or 3\or 4\or 5\or 6\or 7\or 8\or 9\or A\or B\or C\or D\or E\or F\fi}
\mathchardef\mathi="7\hexa\mycmmifam7B
\mathchardef\mathj="7\hexa\mycmmifam7C


\font\mymsbmten=msbm10 at 8pt
\font\mymsbmseven=msbm7 at 5.6pt
\font\mymsbmfive=msbm5 at 4pt
\newfam\mymsbmfam
\textfont\mymsbmfam=\mymsbmten
\scriptfont\mymsbmfam=\mymsbmseven
\scriptscriptfont\mymsbmfam=\mymsbmfive

\mathchardef\mybeth="7\hexa\mymsbmfam69
\mathchardef\mygimmel="7\hexa\mymsbmfam6A
\mathchardef\mydaleth="7\hexa\mymsbmfam6B


\def\placelabel[#1][#2]#3{{%
\setbox10=\hbox{\raise #2cm \hbox{\hskip #1cm #3}}%
\ht10=0pt%
\dp10=0pt%
\wd10=0pt%
\box10}}%


\newif\ifinproclaim%
\global\inproclaimfalse%
\def\proclaim#1{%
\medskip%
%
%
\bgroup%
\inproclaimtrue%
\setbox10=\vbox\bgroup\leftskip=0.8em\noindent{\bf #1}\sl%
}

\def\endproclaim{%
\egroup%
\setbox11=\vtop{\noindent\vrule height \ht10 depth \dp10 width 0.1em}%
\wd11=0pt%
\setbox12=\hbox{\copy11\kern 0.3em\copy11\kern 0.3em}%
\wd12=0pt%
\setbox13=\hbox{\noindent\box12\box10}%
\noindent\unhbox13%
\egroup%
\medskip\ignorespaces%
}

\def\proclaim#1{%
\medskip%
\bgroup%
\inproclaimtrue%
\noindent{\bf #1}%
\nobreak\medskip%
\sl%
}

\def\endproclaim{%
\mypar\egroup\penalty\proclaimpenalty\medskip\ignorespaces%
}

\def\noskipproclaim#1{%
\medskip%
\bgroup%
\inproclaimtrue%
\noindent{\bf #1}\nobreak\sl%
}

\def\endnoskipproclaim{%
\mypar\egroup\penalty\proclaimpenalty\medskip\ignorespaces%
}


\def\proof{{\noindent\bf Proof:\ }}

\def\minf{\mathop{{\mathroman Inf}}}

\def\psl{\Bbb{P}SL}
\def\qed{~$\square$}
\def\munion{\mathop{\cup}}
\def\minter{\mathop{\cap}}
\def\myitem#1{%
    \noindent\hbox to .5cm{\hfill#1\hss}
}

\catcode`\@=11
\def\Eqalign#1{\null\,\vcenter{\openup\jot\m@th\ialign{%
\strut\hfil$\displaystyle{##}$&$\displaystyle{{}##}$\hfil%
&&\quad\strut\hfil$\displaystyle{##}$&$\displaystyle{{}##}$%
\hfil\crcr #1\crcr}}\,}
\catcode`\@=12

\def\makeop#1{%
\global\expandafter\def\csname op#1\endcsname{{\mathroman #1}}}%

\def\makeopsmall#1{%
\global\expandafter\def\csname op#1\endcsname{{\mathroman{\lowercase{#1}}}}}%

\makeopsmall{ArcTan}%
\makeopsmall{ArcCos}%
\makeop{Arg}%
\makeop{Det}%
\makeop{Log}%
\makeop{Re}%
\makeop{Im}%
\makeop{Dim}%
\makeopsmall{Tan}%
\makeop{Ker}%
\makeopsmall{Cos}%
\makeopsmall{Sin}%
\makeop{Exp}%
\makeopsmall{Tanh}%
\makeop{Tr}%
\makeop{End}%
\makeop{Long}%
\makeop{Ch}%
\makeop{Exp}%
\makeop{Int}%
\makeop{Ext}%
\makeop{Aire}%
\makeop{Im}%
\makeop{Conf}%
\makeop{Exp}%
\makeop{Mod}%
\makeop{Log}%
\makeop{Ext}%
\makeop{Int}%
\makeop{Dist}%
\makeop{Aut}%
\makeop{Id}%
\makeop{SO}%
\makeop{Homeo}%
\makeop{Vol}%
\makeop{Ric}%
\makeop{Hess}%
\makeop{Euc}%
\makeop{Isom}%
\makeop{Max}%
\makeop{Long}%
\makeop{Fixe}%
\makeop{Wind}%
\makeop{Mush}%
\makeop{Ad}%
\makeop{loc}%
\makeop{Len}%
\makeop{Area}%
\makeop{SL}%
\makeop{GL}%
\makeop{dVol}%
\makeop{Min}%
\makeop{Symm}%
\makeop{O}%
\makeopsmall{CoTanh}
\makeopsmall{ArcTanh}
\makeopsmall{Sinh}
\makeopsmall{Cosh}
\makeopsmall{ArcCoTanh}

\let\emph=\bf

\hyphenation{quasi-con-formal}

%
%

\ifmakebiblio%
    \openout\biblio=biblio.tex %
    {%
        \edef\fileoutput{\write\biblio{\bgroup\leftskip=2em}}%
        \fileoutput
    }%
\fi%

\newref{LabB}{Labourie F., Probl\`eme de Minkowski et surfaces \`a courbure constante dans les vari\'et\'es hyperboliques, {\sl Bull. Soc. Math. France} {\bf 119} (1991), no. 3, 307--325}
\newref{LabA}{Labourie F., Un lemme de Morse pour les surfaces convexes, {\sl Invent. Math.} {\bf 141} (2000), 239--297}
\newref{RosSpruck}{Rosenberg H., Spruck J. On the existence of convex hyperspheres of constant Gauss curvature in hyperbolic space, {\sl J. Diff. Geom.} {\bf 40} (1994), no. 2, 379--409}
\newref{SmiB}{Smith G., Hyperbolic Plateau problems, Preprint, Orsay (2005)}%
\newref{SmiC}{Smith G., Pointed k-surfaces, to appear in {\it Bull. Soc. Math. France.\/}, math.DG/0506233}%
\newref{SmiG}{Smith G., Finite Area and Volume of Pointed k-Surfaces, in preparation}
\medskip

\ifmakebiblio%
    {\edef\fileoutput{\write\biblio{\egroup}}%
    \fileoutput}%
\fi%

%
%
%
\document
\myfontdefault
\global\chapno=1
\global\showpagenumflag=1
\def\Pagetitle{}
\null
\vfill
\def\centre{\rightskip=0pt plus 1fil \leftskip=0pt plus 1fil \spaceskip=.3333em \xspaceskip=.5em \parfillskip=0em \parindent=0em}%
\def\textmonth#1{\ifcase#1\or January\or Febuary\or March\or April\or May\or June\or July\or August\or September\or October\or November\or December\fi}
\font\abstracttitlefont=cmr10 at 14pt
{\abstracttitlefont\centre A Brief Note on Foliations of Constant Gaussian Curvature\par}
\bigskip
{\centre Graham Smith\par}
\bigskip
{\centre \the\day\ \textmonth\month\ \the\year\par}
\bigskip
{\centre Max Planck Institute for Mathematics in the Sciences,\par
Inselstrasse 22.,\par
D-04103 Leipzig,\par
GERMANY\par}
\bigskip
\noindent{\emph Abstract:\ }This note provides an alternative proof of a result of Labourie. We show that the two complements of the convex 
core of a three dimensional quasi-fuchsian hyperbolic manifold may be foliated by embedded hypersurfaces of constant Gaussian curvature.
\par
\vfill
\nextoddpage
\def\Pagetitle{\sl Foliations of Constant Gaussian Curvature}
%
\newhead{Introduction}
\noindent Let $M$ be a geometrically finite $3$ dimensional hyperbolic manifold without cusps. In
\cite{LabB}, Labourie proved that the complement of the convex core of such a manifold may be foliated by surfaces of constant Gaussian curvature $k$ taking values in the interval $]0,1[$. This permits the
construction of a smooth family of submersions parametrised by $k$ from the moduli space of $M$ into the
Cartesian product of finitely many Teichm\"uller spaces. In this letter, we present a proof of this result in the simpler case where $M$ is a quasi-Fuchsian manifold which may be easily be generalised to higher dimensions.
\medskip
\noindent Let $\Bbb{H}^3$ be $3$ dimensional hyperbolic space. We recall that the ideal boundary of $\Bbb{H}^3$ may be identified with the Riemann sphere. Let $\Gamma$ be a 
cocompact Fuchsian group and let $\rho:\Gamma\rightarrow\opIsom(\Bbb{H}^3)$ be a quasi-fuchsian representation. The quotient 
space $M=\Bbb{H}^3/\rho(\Gamma)$ is a non compact hyperbolic manifold. The limit set of $\rho(\Gamma)$ in 
$\partial_\infty\Bbb{H}^3=\hat{\Bbb{C}}$ is a Jordan curve $c$ which is invariant under the action of $\rho(\Gamma)$. 
Let $C$ be the convex hull of $c$ in $\Bbb{H}^3$. This is also invariant under the action of $\rho{\Gamma}$ and thus 
quotients down to a compact convex subset of $M$, which we will refer to as the convex core of $M$ and will also denote 
by $C$. The complement of the convex core consists of two open non-compact connected components.
\medskip
\noindent The complement in $\hat{\Bbb{C}}$ of the curve $c$ consists of two invariant simply connected domains. Let 
$\Omega$ be one of these domains. $\Omega$ is the interior of the ideal boundary of one of the connected components 
of the complement of $C$ in $\Bbb{H}^3$. Let us denote this component by $U$. $\Omega$ defines a Plateau problem in the
sense of Labourie \cite{LabA}. Thus, using results of Labourie, \cite{LabA}, or Rosenberg and Spruck, \cite{RosSpruck}, for all $k\in]0,1[$, there exists a unique convex embedded submanifold
$\Sigma_k = (S,i_k)$  of constant Guassian curvature equal to $k$ such that if $\hat{\mathi}_k$ is the Gauss lifting
of $i_k$ and if $\overrightarrow{n}:U\Bbb{H}^3\rightarrow\hat{\Bbb{C}}$ is the Gauss Minkowski mapping, then 
$\varphi_k = \overrightarrow{n}\circ\hat{\mathi}_k : S\rightarrow \Omega$ is a homeomorphism. By \cite{SmiB},
the immersion $i_k$ varies continuously with $k$.
Moreover, the ideal boundary of $\Sigma_k$ in $\partial_\infty\Bbb{H}^3$ coincides with $c$. 
\medskip
\noindent We obtain the following result:
\proclaim{Theorem \nextprocno}
\noindent $(\Sigma_k)_{k\in]0,1[}$ defines a smooth foliation of $U$ invariant under the action of $\rho{\Gamma}$.
Moreover, $\Sigma_k$ converges towards $\Omega$ and $\partial U = \overline{U}\minter\overline{C}$ in the Haussdorff 
topology in $\Bbb{H}^3\munion\partial_\infty\Bbb{H}^3$ as $k$ tends towards $1$ and $0$ respectively.
\proclabel{Smooth_Foliation}
\endproclaim
\proof This result follows directly from the lemmata that we will establish in the following sections. Firstly,
by Lemma \procref{It_Is_A_Foliation_Of_Something}, $(\Sigma_k)_{k\in]0,1[}$ is a continuous foliation of some subset of $U$.
Lemmata \procref{The_Foliation_Is_Complete_At_Infinity} and \procref{The_Foliation_Is_Complete_At_Zero} allow us to
establish the convergence properties of this foliation as $k$ tends towards $1$ and $0$ respectively. This allows us in turn
to show that the foliation covers the whole of $U$. Finally, by Lemma \procref{The_Foliation_Is_Smooth}, the foliation is
smooth.\qed
\medskip
\noindent By uniqueness, for all $k$, $\Sigma_k$ is invariant under the action of $\rho(\Gamma)$. Consequently, this foliation projects down to a smooth foliation of one of the connected components of the 
complement of the convex core of $M=\Bbb{H}^3/\rho(\Gamma)$, having analogous convergence properties as $k$
tends to $1$ and $0$.
\medskip 
\noindent Since the result of Rosenberg and Spruck is the only active ingredient of the proof, this result also holds in higher dimensions.
\medskip
\noindent The study of the geometry of this foliation presents many interesting questions. Firstly, we would like to know how the distance from the surface $\Sigma_k$ to the convex core varies with $k$. Indeed, 
by the proof of Lemma \procref{The_Foliation_Is_Complete_At_Zero}, we see that this distance is bounded above
by $\opArcTanh(\sqrt{k})$. Moreover, by considering the Fuchsian case, where 
$\rho{\Gamma}\subseteq\psl(2,\Bbb{R})$ and the convex core is a totally geodesic submanifold, we see that
this bound is strict. However, following the proof of Lemma \procref{The_Foliation_Is_Complete_At_Infinity},
a bound from below is more subtle. One would anticipate, however, that such a bound may be obtained as a 
function of the complex dilatation or Schwarzian derivative of $\Omega$, or of the Haussdorf dimension of $c$.
\medskip
\noindent Secondly, we would like to study the behaviour of the area of $\Sigma_k$ and the volume contained between $\Sigma_k$ and the convex core as $k$ tends towards each of the limits. This corresponds to the study of the asymptotique behaviour of the integrals of $1$ and $f$ over $\Sigma_k$ as $k$ tends to $0$ and $1$. In particular, the constant terms in the asymptotic expansions of these functions would yield a renormalised area and volume.
\medskip
\noindent Thirdly, we observe that an analogous construction may be carried out in the case of pointed $k$-surfaces, as in \cite{SmiC}, although in this case we obtain a lamination rather than a foliation. In particular, in \cite{SmiG}, we show that the area and the volume of pointed $k$-surfaces are well defined and finite. Thus analogous geometrical questions may also be studied in this case.
\medskip
\noindent As a final observation, the immersed hypersurface $\Sigma_k$ carries two different canonical
metrics, being the pull back through $i$ of the hyperbolic metric of $\Bbb{H}^3$, the pull back through 
$\hat{\mathi}$ of the metric of $U\Bbb{H}^3$. It also carries a canonical conformal structure, being the 
pull back through $\overrightarrow{n}\circ\hat{\mathi}$ of the conformal structure of $\hat{\Bbb{C}}$. Using, for example, results proven in \cite{SmiC}, one may show that the Teichm\"uller distance separating these three conformal structures converges to zero as $k$.
\medskip
\noindent This paper is arranged as follows. In the following section, we calculate the effect of a normal deformation on the Gaussian curvature of a surface, and in the remaining sections, we prove the different parts of Theorem \procref{Smooth_Foliation}.
\newhead{Normal Deformation of Immersed Surfaces}
\noindent Let $M$ be a three dimensional hyperbolic manifold and let $\Sigma = (S,\varphi)$ be an immersed hypersurface. We begin by
calculating the effect on the second fundamental form of $\Sigma$ of a variation in the normal direction. 
\medskip
\noindent Let $\Phi:S\times]-\epsilon,\epsilon[\rightarrow M$ be an isotopy of immersions. That is, for all $t\in]-\epsilon,\epsilon[$,
$\Phi_t=\Phi(\cdot,t)$ is an immersion. Suppose that $\Phi_0=\varphi$. Suppose moreover, that the variation is normal. That
is, $\partial_t$ is orthogonal to the tangent space to $S$ with respect to the induced metric obtained by pulling back the
hyperbolic metric on $M$. Let ${\mathsf N}$ be the vector field over $\Phi$ such that, for all $t$, 
${\mathsf N}_t = {\mathsf N}(\cdot,t)$ is the exterior unit normal vector field to $\Sigma_t = (S,\Phi_t)$. 
Let $f_t:S\rightarrow\Bbb{R}$ be such that:
$$
\Phi_*\partial_t = f_t {\mathsf N}_t.
$$
\noindent Finally, for all $t$, let $A_t$ be the Weingarten operator of $\Sigma_t$. Thus, for all $X$ in $TS$:
$$
A\cdot X = \nabla_X{\mathsf N}.
$$
\noindent We obtain the following result:
\proclaim{Lemma \nextprocno}
\noindent The first derivative of $A$ is given by:
\proclabel{First_Derivative_Of_A}
$$
\partial_t A = f\opId - \opHess(f) - fA^2.
$$
\endproclaim
\proof Since this is a local formula, we may assume that $\Phi$ is a diffeomorphism. We thus identify $M$ with
$S\times]-\epsilon,\epsilon[$ and $\Phi$ with the identity. 
\medskip
\noindent The first step involves calculating $\nabla_{\partial_t}{\mathsf N}$. Since the norm of ${\mathsf N}$ is 
constant, we have:
$$
\langle \nabla_{\partial_t} {\mathsf N},{\mathsf N}\rangle = 
1/2\partial_t\langle {\mathsf N},{\mathsf N}\rangle = 0.
$$
\noindent Likewise, if $X$ is tangent to $S$, $\langle {\mathsf N}, X\rangle$ always vanishes, and thus:
$$
\langle \nabla_{\partial_t} {\mathsf N},X\rangle 
= -\langle {\mathsf N},\nabla_{\partial_t} X\rangle
= -\langle {\mathsf N},\nabla_X f{\mathsf N}\rangle.
$$
\noindent Since $\nabla_X{\mathsf N}$ is tangent to $S$, this yields:
$$
\langle \nabla_{\partial_t} {\mathsf N},X\rangle 
= -\langle \nabla f, X\rangle.
$$
\noindent It thus follows that:
$$
\nabla_{\partial_t} {\mathsf N} = -\nabla^\Sigma f.
$$
\noindent Where $\nabla^\Sigma f$ is the projection onto $T\Sigma$ of the gradiant of $f$. Now, for all $X,T\in TS$, we have:
$$\matrix
\partial_t\langle \nabla_X{\mathsf N}, Y\rangle\hfill
&=\langle\nabla_{\partial_t}\nabla_X{\mathsf N},Y\rangle + \langle\nabla_X{\mathsf N},\nabla_{\partial_t}Y\rangle\hfill\cr
&=\langle R_{{\partial_t}X}{\mathsf N},Y\rangle + \langle \nabla_X\nabla_{\partial_t}{\mathsf N},Y\rangle
+\langle \nabla_X{\mathsf N}, \nabla_Y f{\mathsf N}\rangle.\hfill\cr
&=\langle fR_{{\mathsf N}X}{\mathsf N},Y\rangle - \langle \nabla_X\nabla^\Sigma f, Y\rangle + \langle A\cdot X, A\cdot Y\rangle.\hfill\cr
\endmatrix$$
\noindent Observe that the above formula contains not only the variation of $A$ but also the variation of the metric. We thus
have to subtract the variation of the metric in order to obtain the variation of $A$. For each $t$, we denote by $g_t$ the
restriction of $\langle\cdot,\cdot\rangle$ to the surface $S_t = S\times\{t\}$. We have:
$$\matrix
\partial_tg_t(X,Y) \hfill&=g_t(\nabla_{\partial_t}X,Y) + g_t(X, \nabla_{\partial_t}Y)\hfill\cr
&=g_t(\nabla_X {\partial_t},Y) + g_t(X, \nabla_Y {\partial_t})\hfill\cr
&=f g_t(A\cdot X, Y) + f g_t(X, A\cdot Y).\hfill\cr
\endmatrix$$
\noindent Combining this with the preceeding formula, and bearing in mind that $A$ is symmetric, we obtain:
$$
g_t((\partial_t A_t)\cdot X, Y) = g_t(fR_{{\mathsf N}X}{\mathsf N} - \opHess(f)\cdot X - fA^2,Y).
$$
\noindent Finally, since $M$ is hyperbolic, for all $X$ tangent to $S$, the curvature satisfies:
$$
R_{{\mathsf N}X}{\mathsf N} = X.
$$
\noindent The result now follows.\qed
\medskip
\noindent In particular, when $f\cong 1$, we obtain the family of immersed surfaces equidistant from 
$\Sigma$. Let us note by $\Sigma_t = (S,i_t)$ this family, and for all $t$ let $A_t$ be the Weingarten operator of 
$\Sigma_t$. We obtain the following corollary:
\proclaim{Corollary \nextprocno}
\noindent $A_t$ satisfies the following differential equation:
\proclabel{Diff_Eq_For_Equidistant_Family}
$$
\partial_t A_t = \opId - A_t^2.
$$
\endproclaim
\noindent This permits us to calculate $A_t$ over all time. Indeed, let $\lambda_1(t)\leqslant\lambda_2(t)$ be the 
two eigenvalues of $A_t$. We have the following result:
\proclaim{Corollary \nextprocno}
\noindent For each $i$, there exists a constant $T_0\in\Bbb{R}$ such that, either $\lambda_i(t)=\opTanh(T_0+t)$ or
$\lambda_i(t)=\opCoTanh(T_0+t)$ or $\lambda_i(t)=1$.
\proclabel{Growth_Is_Assured}
\endproclaim
\proof Let $p$ be a point in $S$. Let $E_1,E_2$ be the principal eigenvectors of $A_0$ at $p$. If we now define the
symmetric matrix $\hat{A}_t$ such that $E_1$ and $E_2$ are its eigenvectors and its eigenvalues are given by the
formulae in the statement of this corollary, then $\hat{A}_t$ satisfies the same differential equation as $A_t$ with the
same initial conditions, and they thus coincide. The result now follows.\qed
\medskip
\noindent This permits us to obtain the following bound from above:
\proclaim{Lemma \nextprocno}
\noindent Let $k\in ]0,1[$ and suppose that $\Sigma$ is convex and of constant Gaussian curvature equal to $k$. Suppose,
moreover, that there exists $a>0$ such that, in the sense of positive definite matrices:
$$
A_0 \geqslant a\opId.
$$
\noindent Then, there exists a function $K:[0,+\infty[\rightarrow[0,1]$ which only depends on $k$ and $a$ such that:
\proclabel{Curvature_Of_Equidistant_Surfaces_Bounded_Above}
\medskip
\myitem{(i)} $K$ is strictly increasing,
\medskip
\myitem{(ii)} $K(t)$ tends to $1$ as $t$ tends to $+\infty$, and
\medskip
\myitem{(iii)} for all $t$, the Gaussian curvature of $\Sigma_t$ is bounded above by $K(t)$.
\endproclaim
\proof Let $p$ be a point in $S$. For all $t$, let $\lambda_1(t)\leqslant\lambda_2(t)$ be the eigenvalues of
$A_t(p)$. Since $\lambda_1>0$ and $0<\lambda_1\lambda_2<1$, by Corollary \procref{Growth_Is_Assured}, there exists 
$A_1 > \hat{a} := \opArcTanh(a)$ such that:
$$
\lambda_1(t) = \opTanh(A_1 + t).
$$
\noindent Corollary \procref{Growth_Is_Assured} now offers three possibilities for the evolution of $\lambda_2$, depending
on the value of $A_1$:
\medskip
\myitem{(i)} If $A_1<\opArcTanh(k)$ then $\lambda_1(0)<k$ and so $\lambda_2(0)>1$. Thus, there exists $A_2>0$ such that:
$$
\lambda_2(t) = \opCoTanh(A_2 + t).
$$
\noindent Since $\lambda_1\lambda_2<1$, $A_2>A_1$. Consider the function $\varphi(t) = \opCoTanh(c+t)\opTanh(t)$ where
$c>0$. Using elementary hyperbolic trigonometry, we obtain:
$$
\partial_t\varphi(c,t) = \opSinh(2c)/2(\opCosh(t)\opSinh(c+t))^2 > 0.
$$
\noindent Consequently, $\varphi(c,t)$ is increasing in $t$. Moreover, $\varphi(c,t)$ tends to $1$ as $t$ tends to 
infinity. Since $\varphi((A_2-A_1) + A_1,A_1) = k < 1$ and $A_1>\hat{a}$:
$$\matrix
&\varphi(A_2 - A_1 + \hat{a}, \hat{a}) \hfill&\leqslant k\hfill\cr
\Rightarrow\hfill& A_2 - A_1 \hfill&\geqslant\opArcCoTanh(ka) - \hat{a}\hfill\cr
\Rightarrow\hfill& A_2 \hfill&\geqslant\opArcCoTanh(ka).\hfill\cr
\endmatrix$$
\noindent We now define $\hat{A}_1(a,k) = \opArcTan(k)$ and $\hat{A}_2(a,k) = \opArcCoTanh(ka)$. We define 
$K_1(a,k;t):[0,\infty[\rightarrow\Bbb{R}$ by:
$$\matrix
K_1(a,k;t) \hfill&= \opTanh(\hat{A}_1(a,k) + t)\opCoTanh(\hat{A}_2(a,k) + t)\hfill\cr
& = \varphi(\hat{A}_2(a,k) - \hat{A}_1(a,k), t + \hat{A}_1(a,k)).\hfill\cr
\endmatrix$$
\noindent The function $K_1$ is increasing in $t$ and tends to $1$ as $t$ tends to infinity. Moreover, since 
$\opTanh$ is increasing and $\opCoTanh$ is decreasing, bearing in mind the bounds on $A_1$ and $A_2$, for all $t$:
$$
K_1(a,k;t) \geqslant \lambda_1(t)\lambda_2(t).
$$
\noindent $K_1$ is thus the desired function in this case.
\medskip
\myitem{(ii)} If $A_1=\opArcTanh(k)$, then $\lambda_1(0)=k$ and so $\lambda_2(0)=1$. Thus, for all $t$:
$$
\lambda_2(t) = 1.
$$
\noindent We define $K_2(a,k;t) = \opTanh(\opArcTanh(k) + t)$, and we see that $K_2$ is the desired function in this case.
\medskip
\myitem{(iii)} If $A_2>\opArcTanh(k)$, then $\lambda_1(0)>k$ and so $\lambda_2(0)<1$. Thus, there exists $A_2>A_1$ such that:
$$
\lambda_2(t) = \opTanh(A_2 + t).
$$
\noindent We define $\hat{A}_1=\opArcTanh(\sqrt{k})$ and $\hat{A}_2=\opArcTanh(k/a)$. We have:
$$
A_1 \leqslant \hat{A}_1,\qquad A_2 \leqslant \hat{A}_2.
$$
\noindent We define the function $K_3(a,k;t)$ by:
$$
K_3(a,k;t) = \opTanh(\hat{A}_1 + t)\opTanh(\hat{A}_2 + 2).
$$
\noindent We see that $K_3$ is the desired function in this case. The result now follows by taking the maximum of 
$K_1$, $K_2$ and $K_3$.\qed
\newhead{The Family is a Foliation}
\proclaim{Lemma \nextprocno}
\noindent For $k\neq k'$, $\Sigma_k$ and $\Sigma_{k'}$ are disjoint.
\proclabel{It_Is_A_Foliation_Of_Something}
\endproclaim
\proof We may assume that $k>k'$. For $t\in [0,1]$ let $\Omega_t$ be a family of Jordan subdomains of $\Omega$ such that:
\medskip
\myitem{(i)} $\Omega_0$ is a disc,
\medskip
\myitem{(ii)} for $t<t'$ the closure of $\Omega_t$ is contained in $\Omega_{t'}$, and
\medskip
\myitem{(iii)} $\partial\Omega_t$ tends to $c=\partial\Omega$ as $t$ tends to $1$.
\medskip
\noindent Such a family may be obtained, for example, by uniformising the annulus obtained by removing a disc from $\Omega$. 
Let $\delta$ be a real number such that $0<\delta<1-k$. For $t\in[0,1]$ we define $k_t\in [k,1]$ by:
$$
k_t = (1-\delta)(1-t) + kt.
$$
\noindent For all $t$, let $\Sigma'_t = (S,j_t)$ be the unique solution to the Plateau problem given by $\Omega_t$ with
constant Gaussian curvature equal to $k_t$. We may assume that we have have chosen $\Omega_0$ sufficiently small so that 
$\Sigma'_0$ is disjoint from $\Sigma_{k'}$ in $\Bbb{H}^3\munion\partial_\infty\Bbb{H}^3$. 
\medskip
\noindent We define $t_0\in[0,1]$ by:
$$
t_0 = \minf\{t\in[0,1]\text{ s.t. } \Sigma'_t\minter\Sigma_{k'}\neq\emptyset\text{ in }\Bbb{H}^3\munion\partial_\infty\Bbb{H}^3.\}.
$$
\noindent Suppose that this set is non-empty, and thus that $t_0\in[0,1]$. By compactness, $\Sigma'_{t_0}$ and 
$\Sigma_{k'}$ intersect non-trivially in $\Bbb{H}^3\munion\partial_\infty\Bbb{H}^3$. Since the ideal boundaries of 
$\Sigma'_{t_0}$ and $\Sigma_{k'}$ are $\partial\Omega_{t_0}$ and $\partial\Omega$ respectively, and since these do not
intersect, it follows that these two surfaces intersect non-trivially in $\Bbb{H}^3$. However, by continuity, since
$\Sigma_0$ lies in (the closure of) the exterior of $\Sigma_{k'}$, so does $\Sigma'_{t_0}$. Consequently, $\Sigma'_{t_0}$
is an exterior tangent to $\Sigma_{k'}$. This is absurd by the geometric maximum principal, 
since the Gaussian curvature
of $\Sigma'_{t_0}$ exceeds that of $\Sigma_{k'}$. Consequently, the set used to define $t_0$ is empty, and the result now
follows.\qed
\newhead{Completeness of the Foliation}
\noindent We recall that $\Bbb{H}^3\munion\partial_\infty\Bbb{H}^3$ has the topology of a closed ball. We now prove 
completeness of the foliation near infinity:
\proclaim{Lemma \nextprocno}
\noindent $(\Sigma_k)_{k\in[0,1]}$ tends towards $\Omega$ in the Haussdorff topology as $k$ tends to infinity.
\proclabel{The_Foliation_Is_Complete_At_Infinity}
\endproclaim
\proof For all $k$, let $A_k$ be the Weingarten operator of $\Sigma_k$. Let $k_0\in[0,1]$ be arbitrary. Since 
$\rho(\Gamma)$ has a cocompact action of $\Sigma_{k_0}$, there exists $a>0$ such that, in the sense of positive definite
symmetric matrices, for all $p\in S$:
$$
A_k(p) \geqslant a\opId.
$$
\noindent Let the function $K(a,k_0;\cdot):[0,+\infty[\rightarrow\Bbb{R}$ be as in 
Lemma \procref{Curvature_Of_Equidistant_Surfaces_Bounded_Above}. For all $d\in [0,\infty[$, let $\Sigma'_d = (S,j_d)$ be
equidistant surface of distance $d$ from $\Sigma_{k_0}$. For all $d$, let $k'_d:S\rightarrow ]0,\infty[$ be the Gaussian
curvature of $\Sigma'_d$. By Lemma \procref{Curvature_Of_Equidistant_Surfaces_Bounded_Above}, for all $d$, and for
all $p\in S$:
$$
k'_d(p) \leqslant K(a,k_0;d) < 1.
$$
\noindent Let $D>0$ be arbitrary. By a continuity argument identical to that employed in the proof of Lemma
\procref{It_Is_A_Foliation_Of_Something}, for all $k\geqslant K(a,k_0;D)$, the immersed surface
$\Sigma_k$ lies in the exterior of $\Sigma'_D$. Since $D$ may be chosen arbitrarily large, the result now follows.\qed
\medskip
\noindent In order to prove the completeness of the foliation near the convex core, we require the following definition
of Gaussian curvature in the weak sense:
\proclaim{Definition \nextprocno}
\noindent Let $(N,i)$ be a (not necessarily $C^2$) hypersurface in a manifold $M$. Let $p$ be a point in $N$. We say that
$(N,i)$ is weakly convex at $p$ if there exist open neighbourhoods $U\subseteq N$, $V\subseteq M$ of $p$ and $i(p)$
respectively and a convex subset $K\subseteq V$ such that $i(p)\in\partial K$ and, for all $q\in U$, $i(q)\in K$.
\medskip
\noindent If $(N,i)$ is weakly convex at $p$, for $k>0$, we say that the Gaussian curvature of $(N,i)$ is at least $k$ in the
weak sense if and only if $U$, $V$ and $K$ may be chosen such that $\partial K$ is smooth and that the Gaussian curvature of
$K$ at $i(p)$ is equal to $k$.
\endproclaim
\noindent For all $d$, let $\Sigma'_d = (S,j_d)$ be the equidistant surface in $U$ of distance $d$ from the convex core. We
have the following result:
\proclaim{Lemma \nextprocno}
\noindent For all $d\in [0,+\infty[$ and for all $p\in S$, $\Sigma'_d$ is weakly convex at $p$ and its Gaussian curvature
is at least $\opTanh(d)^2$ in the weak sense.
\proclabel{Weak_Lower_Bound_Of_Curv_For_Level_Surfaces}
\endproclaim
\proof Let $d$ be a non negative real number. Let $p$ be a point in $S$. Let $q$ be the point in $\partial U$ closest to
$i_d(p)$. There exists a supporting hyperplane to $C$ at $q$ whose normal points towards $i_d(p)$. Let $P$ be this 
supporting hyperplane. Let $P_d$ be the equidistant surface in $U$ at distance $d$ from $P$. Since $C$ lies entirely to 
one side of $P$, $\Sigma'_d$ lies within the interior of $P_d$. Moreover, since the normal of $P$ at $q$ points towards
$i_d(p)$, $i_d(p)$ lies on $P_d$. However, $P_d$ is strictly convex, and, by Corollary \procref{Growth_Is_Assured}, its
Gaussian curvature equals $\opTanh(d)^2$. The result now follows.\qed
\medskip
\noindent We are now in a position to prove the completeness of the foliation near the convex core:
\proclaim{Lemma \nextprocno}
\noindent $(\Sigma_k)_{k\in[0,1]}$ tends towards $\partial U$ in the Haussdorff topology as $k$ tends to $0$.
\proclabel{The_Foliation_Is_Complete_At_Zero}
\endproclaim
\proof Let $\Omega'$ be the other connected component of $c$ in $\hat{C}$. For $t\in[0,1]$, let $\Omega'_t$ be a family 
of Jordan subdomains of $\Omega'$ such that:
\medskip
\myitem{(i)} $\Omega'_0$ is a disc,
\medskip
\myitem{(ii)} for $t<t'$ the closure of $\Omega'_t$ is contained in $\Omega'_{t'}$, and
\medskip
\myitem{(iii)} $\partial\Omega'_t$ tends to $c=\partial\Omega'$ as $t$ tends to $1$.
\medskip
\noindent Such a family may be obtained, for example, by uniformising the annulus obtained by removing a disc from 
$\Omega'$. Let $\delta$ be a positive real number. Suppose that $k<\delta$. Let $d$ be such that the Gaussian curvature 
of $\Sigma'_d$ is greater than $\delta$ in the weak sense. For $t\in[0,1]$ we define $k_1\in[0,k]$ by:
$$
k_t = t\delta.
$$
\noindent For all $t$, let $\Sigma''_t = (S,k_t)$ be the unique solution to the Plateau problem given by
$\hat{\Bbb{C}}\setminus\overline{\Omega}_t$ with constant Guassian curvature equal to $k_t$. $\Sigma''_0$ is a totally
geodesic submanifold of $\Bbb{H}^3$. We may assume that $\Omega_0$ is chosen sufficiently small so that $\Sigma''_0$
is disjoint from $\Sigma'_d$ in $\Bbb{H}^3\munion\partial_\infty\Bbb{H}^3$. By completing the continuity argument as in the
proof of Lemma \procref{It_Is_A_Foliation_Of_Something}, we see that $\Sigma_k$ lies in the interior of $\Sigma'_d$. Since
$d$ may be made arbitarily small, the result now follows.\qed
\newhead{Smoothness of the Foliation}
\noindent We have the following result:
\proclaim{Lemma \nextprocno}
\noindent The foliation $(\Sigma_k)_{k\in]0,1[}$ is smooth.
\proclabel{The_Foliation_Is_Smooth}
\endproclaim
\proof We consider the quotient foliation in $\Bbb{H}^3/\rho(\Gamma)$. Thus every surface in the foliation is compact.
Let $k$ be a real number in $]0,1[$. Let $M$ be a field of symmetric positive definite matrices over $S$. For $f$
a function over $S$, we define:
$$
\Delta^M f = \opTr(M^{-1}\opHess(f)).
$$
\noindent Consider now the following system of coupled partial differential equations:
$$\matrix
\partial_t A \hfill&= f\opId - \opHess(f) - fA^2,\hfill\cr
\Delta^Af + \opTr(A - A^{-1})f \hfill&= -1,\hfill\cr
\endmatrix$$ 
\noindent subject to the initial condition that $A(0)=A_k$. The second condition ensures that $\opDet(A(t)) = k + t$ for all $f$. Let 
$M$ be a $2\times 2$ symmetric matrix such that $\opDet(M)=k$. If $\lambda\in\Bbb{R}$ be the lower of the two eigenvalues
of $M$, then:
$$
\opTr(A - A^{-1}) = \frac{(k-1)(\lambda^2 + k)}{\lambda k}.
$$
\noindent Since $\opDet(A_k) = k < 1$, it follows that for small values of $t$:
$$
\opTr(A_t - A_t^{-1}) < 0.
$$
\noindent Consequently, by the maximum prinicipal, since $S$ is compact, for small $t$, the operator 
$\Delta^{A_t} + \opTr(A_t-A_t^{-1})$ is injective. Any Laplacian acting on the space of functions over a compact manifold 
is of index zero, and thus, for small $t$, the operator $\Delta^{A_t} + \opTr(A_t-A_t^{-1})$ is surjective. Thus, for all
$t$, there exists a unique solution $f_t$ to the second equation. Using classical techniques of partial differential 
equations, we may thus smoothly solve this system for small values of $t$.
\medskip
\noindent We now consider the mapping $\Phi:S\times]-\epsilon,\epsilon[\rightarrow M$ and the vector field
${\mathsf N} \in \Phi^*TM$ defined by the following system of coupled ordinary differential equations:
$$\matrix
\partial\Phi \hfill&= f{\mathsf N},\hfill\cr
\nabla_{\partial_t} {\mathsf N} \hfill&= -\nabla^\Sigma f,\hfill\cr
\endmatrix$$
\noindent subject to the condition that $\Phi(\cdot,0)=i_k$ and that ${\mathsf N}(\cdot,0)$ is the exterior unit normal
vector field to $\Sigma_k$. This system of equations may be smoothly and uniquely solved for small values of $t$. 
Following the proof of Lemma \procref{First_Derivative_Of_A}, we see that, since they both satisfy the same differential equations with the same initial conditions, for all $t$, the field ${\mathsf N}_t$ is equivalent to the
exterior unit normal vector field to $(S,\Phi_t)$. Likewise, by Lemma \procref{First_Derivative_Of_A}, the Weingarten 
operator of $(S,\Phi_t)$ coincides with $A_t$. Consequently, the surface $(S,\Phi_t)$ is of constant Guassian curvature 
equal to $k+t$, and thus, by uniqueness, it coincides with $\Sigma_{k+t}$. 
\medskip
\noindent It follows that the foliation is smooth near $\Sigma_k$, and since $k\in]0,1[$ is arbitrary, the result now
follows.\qed
\goodbreak
\newhead{Bibliography}
{\leftskip = 5ex \parindent = -5ex
\leavevmode\hbox to 4ex{\hfil \cite{LabB}}\hskip 1ex{Labourie F., Probl\`eme de Minkowski et surfaces \`a courbure constante dans les vari\'et\'es hyperboliques, {\sl Bull. Soc. Math. France} {\bf 119} (1991), no. 3, 307--325}
\medskip
\leavevmode\hbox to 4ex{\hfil \cite{LabA}}\hskip 1ex{Labourie F., Un lemme de Morse pour les surfaces convexes, {\sl Invent. Math.} {\bf 141} (2000), 239--297}
\medskip
\leavevmode\hbox to 4ex{\hfil \cite{RosSpruck}}\hskip 1ex{Rosenberg H., Spruck J. On the existence of convex hyperspheres of constant Gauss curvature in hyperbolic space, {\sl J. Diff. Geom.} {\bf 40} (1994), no. 2, 379--409} 
\medskip
\leavevmode\hbox to 4ex{\hfil\cite{SmiB}}\hskip 1ex{Smith G., Hyperbolic Plateau problems, Preprint, Orsay (2005)}%
\medskip
\leavevmode\hbox to 4ex{\hfil \cite{SmiC}}\hskip 1ex{Smith G., Pointed k-surfaces, to appear in {\it Bull. Soc. Math. France.\/}, \hfill\break math.DG/0506233}%
\medskip
\leavevmode\hbox to 4ex{\hfil \cite{SmiG}}\hskip 1ex{Smith G., Finite Area and Volume of Pointed k-Surfaces, in preparation}%
\medskip
}%
\enddocument